\documentclass[12pt]{article}
\usepackage{amsmath,amsthm,amscd,amsfonts,amssymb,graphicx}
\def\demo{{\it Proof. }}
\def\fin{\hfill{$\square$}}
\newtheorem{theorem}{Theorem}[section]
\newtheorem{lemma}{Lemma}[section]
\newtheorem{definition}{Definition}[section]

\newcommand\Ker{\operatorname{ker}}

\newcommand\coh{\operatorname{H}}

\newcommand\Spe{\operatorname{Spec}}

\newcommand\Pic{\operatorname{Pic}}

\newcommand\rk{\operatorname{r}}
\newcommand\rg{\operatorname{rk}}
\newcommand\dg{\operatorname{d}}
\newcommand\Jac{\operatorname{\overline{Jac}}}
\newcommand\jac{\operatorname{Jac}}
\newcommand\Ja{\operatorname{\overline{Jac}}}
\newcommand\ja{\operatorname{Jac}}
\newcommand\so{\operatorname{Supp}}
\newcommand\des{\operatorname{\underset{(\leq)}<}}
\newcommand\gen{\operatorname{g}}

\title{Simpson Jacobians of generalized tree-like curves}
\author{Ana Cristina L\'opez Mart\'{\i}n}

\begin{document}
\maketitle
\begin{abstract}

The compactified Jacobian of any projective curve $X$ is defined
as the Simpson moduli space of torsion free rank one degree $d$
sheaves that are semistable with respect to a fixed polarization
$H$ on $X$. In this paper we give explicitly the structure of this
compactified Simpson Jacobian in the case where $X$ is a
generalized tree-like curve, i.e., a projective, reduced and
connected curve such that the intersection points of its
irreducible components are disconnecting ordinary double points.
We prove that it is isomorphic to the product of the compactified
Jacobians of a certain degree $d_i$ of its components $C_i$, where
the degrees $d_i$ depend on $d$, $H$ and on the particular
structure of the curve. We find also necessary and sufficient
conditions for the existence of stable points which allow us to
study the variation of these Simpson Jacobians as the polarization
$H$ changes.
\end{abstract}
\section{Introduction}

The problem of compactifying the generalized Jacobian of a
singular curve has been studied since Igusa's work \cite{I} around
1950. He constructed a compactification of the Jacobian of a nodal
and irreducible curve $X$ as the limit of the Jacobians of smooth
curves approaching $X$. Igusa also showed that his
compactification does not depend on the considered family of
smooth curves. An intrinsic characterization of the boundary
points of the Igusa's compactification as the torsion free, rank 1
sheaves which are not line bundles is due to Mumford and Mayer.
The complete construction for a family of integral curves over a
noetherian Hensel local ring with residue field separably closed
was carried out by D'Souza \cite{D'So}. One year later, Altman and
Kleiman \cite{AK} gave the construction for a general family of
integral curves.

When the curve $X$ is reducible and nodal, Oda and Seshadri
\cite{OS} produced a family of compactified Jacobians
$\jac_{\phi}$ parameterized by an element $\phi$ of a real vector
space. Seshadri dealt in \cite{Ses} with the general case of a
reduced curve considering sheaves of higher rank as well.

In 1994, Caporaso showed \cite{C} how to compactify the relative
Jacobian over the moduli of stable curves and described the
boundary points of the compactified Jacobian of a stable curve $X$
as invertible sheaves on certain Deligne-Mumford semistable curves
that have $X$ as a stable model. Recently, Pandharipande \cite{P}
has given  another construction with the boundary points now
representing torsion free, rank 1 sheaves and he showed that the
Caporaso's compactification was equivalent to his.

On the other hand, Esteves \cite{E} constructed a compactification
of the relative Jacobian of a family of geometrically reduced and
connected curves and compared it with Seshadri's construction
\cite{Ses} using theta funtions and Alexeev  \cite{A} gave a
description of the Jacobian of certain singular curves in terms of
the orientations on complete subgraphs of the dual graph of the
curve.

Simpson's work \cite{S} on the moduli of pure coherent sheaves on
projective spaces allows us to define in a natural way the
Jacobian of any projective curve $X$ as the space $\jac^d(X)_s$ of
equivalence classes of stable invertible sheaves with degree $d$.
This is precisely the definition we adopt and we also denote by
$\Jac^d(X)$ the space of equivalence classes of semistable, pure
dimension 1, rank 1 sheaves with degree $d$.

In some recent papers about the moduli spaces of stable vector
bundles on elliptic fibrations, for instance \cite{HM}, the
Jacobian in the sense of Simpson of spectral curves appears.
Beauville \cite{B} uses it as well in counting the number of
rational curves on K3 surfaces. This suggests the necessity to
determine the structure of these Simpson Jacobians.

The aim of this work is to describe the Simpson Jacobian of a
generalized tree-like curve, that is, a projective, reduced and
connected curve such that the intersection points of its
irreducible components are disconnecting ordinary double points.

Let $H$ be the fixed polarization on $X$. Let $C_1,\hdots, C_N$
denote the irreducible components of $X$ and let
$P_1,\hdots,P_{N-1}$ be the intersection points of $C_1,\hdots,
C_N$. The first result we use is a lemma of Teixidor \cite{T} that
allows to order $C_i$ and to find subcurves $X_i$ of $X$ that are
also generalized tree-like curves and intersect its complement in
$X$ at just one point $P_i$. Then, we write
$k_{X_i}=\frac{h_{X_i}(b+1)}{h}$ for $i=1,\hdots,N-1$, where $h$
is the degree of $H$, $h_{X_i}$ is the degree of the induced
polarization on $X_i$ and $b$ is the residue class of $d$ minus
the arithmetic genus of $X$ modulo $h$.

The theorem that gives the description of the schemes
$\jac^d(X)_s$ and $\Jac^d(X)$ is the following:

{\bf Theorem \ref{t}.} {\it Let $X=C_1\cup\hdots\cup C_N$, $N\geq
2$, be a generalized tree-like curve with a polarization $H$. If
$k_{X_i}$ is not integer for $i\leq N-1$, then $\jac^d(X)_s$ is
equal to $\prod_{i=1}^{N} \Pic^{d_i^X}(C_i)$ and
$$\jac^d(X)_s\subseteq \Jac^d(X)_s=\Jac^d(X)\simeq \prod_{i=1}^{N}
\Ja^{d_i^X}(C_i)$$ where $d_i^X$ are integer numbers inductively
constructed.

If $k_{X_i}$ is integer for some $i\leq N-1$, then $\jac^d(X)_s$
and $\Jac^d(X)_s$ are empty and $$\Jac^d(X)\simeq
\prod_{i=1}^{N}\Ja^{d_i}(C_i)$$ where $d_i$  are integer numbers
recursively constructed with an algorithm.}

As we will see, these integer numbers, $d_i$ and $d_i^X$, depend
only on the degree of $H$ and on the fixed ordering of the
irreducible components of $X$.

The first we prove is that a pure dimension 1 sheaf of rank 1 and
degree $d$ on $X$ that is locally free at $P_i$, $i=1,\hdots,
N-1$, is stable if and only if $k_{X_i}\notin \mathbb{Z}$ for
$i=1,\hdots,N-1$ and $F$ is obtained by gluing torsion free rank 1
sheaves on $C_i$ whose degrees are precisely the integers $d_i^X$
(see lemmas \ref{lem5} and \ref{lem6} if $F$ is a line bundle,
otherwise the proof is the  same). By the contrary, if $F$ is not
locally free at some intersection point $P_i$, the sheaf cannot be
stable. Finally, we prove that strictly semistable torsion free
rank 1 sheaves exist only when $k_{X_i}\in\mathbb{Z}$ for some
$i=1,\hdots,N-1$ and in this case we construct a Jordan-H\"{o}lder
filtration for it to conclude that its $S$-equivalence class
belongs to $\prod_{i=1}^{N}\Ja^{d_i}(C_i)$.

It is important to note that, although
$\Jac^d(C_i)\simeq\Jac^{d'}(C_i)$ for any integers $d$ and $d'$,
because $C_i$ is an integral curve, to obtain a point in
$\Jac^d(X)$ it is necessary to consider sheaves on $C_i$ of
degrees exactly the integer numbers $d_i$

Theorem \ref{t} allows also to describe a wall-crossing structure
for these Simpson Jacobians, in other words, to analyze the
variation of the moduli spaces $\jac^d(X)_s$ as the polarization
changes. The variation of the polarization is studied by
considering the variation of the numbers $k_{X_i}$,
($i=1,\hdots,N-1$). The space of possible values of such numbers
contains a finite number of walls, the hyperplanes $k_{X_i}=a$
where $a$ is an integer number $1\leq a\leq h_{X_i}$. The
complement of the walls has a finite number of connected
components each of them called a chamber. Therefore, Theorem
\ref{t} means that in the interior of a given chamber the Jacobian
$\jac^d(X)_s$ is not empty and independent of the choice of
numbers $k_{X_i}$. When the numbers $k_{X_i}$ lie in a wall, the
theorem shows that the Simpson Jacobian of $X$ becomes empty.
Finally, if we choose two families of numbers $k_{X_i}$ lying in
the interior of two different chambers, then the structure of the
corresponding Simpson Jacobians is the same but the degrees
$d_i^X$ of the induced line bundles on every component $C_i$
change according to \eqref{des1}.

If $\alpha$ is a real number, we use [$\alpha$] to denote the
greatest integer less than or equal to $\alpha$.

\section{Preliminaries}

Let $X$ be a projective, reduced and connected curve over an
algebraically closed field $\kappa$. Let $C_1,\hdots, C_N$ denote
the irreducible components of $X$. Let $\mathcal{L}$ be an ample
invertible sheaf on $X$, let $H$ be the associated polarization
and let $h$ be the degree of $H$.

Let $F$ be a coherent sheaf on $X$. We say that $F$ is pure of
dimension one or torsion free if  for all nonzero subsheaves
$F'\hookrightarrow F$ the dimension of $\so(F')$ is 1. The rank
and the degree (with respect to $H$) of $F$ are the rational
numbers $\rk_H(F)$ and $\dg_H(F)$ determined by the Hilbert
polynomial
$$P(F,n,H)=\chi(F\otimes\mathcal{O}_X(nH))=h\rk_H(F)n+\dg_H(F)+\rk_H(F)\chi(\mathcal{O}_X).$$
The slope of $F$ is defined by
$$\mu_H(F)=\frac{\dg_H(F)}{\rk_H(F)}$$ The sheaf $F$ is stable
(resp. semistable) with respect to $H$ if $F$ is pure of dimension
one and for any proper subsheaf $F'\hookrightarrow F$ one has
$$\mu_H(F')<\mu_H(F) \ (\text{resp.} \leq)$$

In \cite{S} Simpson defined the rank as the integer number
$h\rk_H(F)$ and the slope as the quotient
$$\frac{\dg_H(F)+\rk_H(F)\chi(\mathcal{O}_X)}{h\rk_H(F)}$$
Stability and semistability considered in terms of Simpson's slope
and in terms of $\mu_H$ are equivalent. We adopt these definitions
of rank and degree of $F$ because they coincide with the classical
ones when the curve is integral.

\noindent {\bf Notation:} We shall use the following notation. For
every proper subcurve $D$ of $X$, we will denote by $F_D$ the
restriction of $F$ to $D$ modulo torsion, that is,
$F_D=(F\otimes\mathcal{O}_D)/\text{torsion}$, $\pi_D$ will be the
surjective morphism $F\to F_D$ and $F^D=\Ker\pi_D$. We shall
denote by $h_D$ the degree of the induced polarization $H_D$ on
$D$. If $d=\dg_H(F)$ then we shall write $d_D=\dg_{H_D}(F_D)$.

We first recall some general properties we will use later.

\begin{lemma}\label{lem1} Let $F$ be a pure dimension one sheaf on $X$
supported on a subcurve $D$ of $X$. Then $F$ is stable (resp.
semistable) with respect to $H_D$ if and only if $F$ is stable
(resp. semistable) with respect to $H$.
\end{lemma}

\begin{proof} It follows from the equality $$P(F,n,H)=\chi(i_*F\otimes
\mathcal{O}_X(nH))=\chi(F\otimes \mathcal{O}_D(nH_D))=P(F, n,
H_D)$$ where $i\colon D\hookrightarrow X$ is the inclusion map.
\end{proof}

\begin{lemma} \label{lem2} A torsion free  rank 1 sheaf $F$ on $X$ is stable (resp.
semistable) if and only if $\mu_H(F^D)<\mu_H(F)$ (resp. $\leq$)
for every proper subcurve $D$ of $X$.
\end{lemma}

\begin{proof} Given a subsheaf $G$ of $F$ such that $\so(G)=D\subset X$,
let us consider the complementary subcurve $\overline{D}$ of $D$
in $X$, i.e. the closure of $X-D$. Since $F_{\overline D}$ is
torsion free, we have $G\subset F^{\overline D}$ with
$\rk_H(G)=\rk_H(F^{\overline D})$ so that
$\mu_H(G)\leq\mu_H(F^{\overline D})$ and the result follows.
\end{proof}

\bigskip
Let $P_1,\hdots, P_k$ denote the intersection points of
$C_1,\hdots,C_N$. It is known (see \cite{Ses}) that for every pure
dimension one sheaf $F$ on $X$ there is an exact sequence $$0\to
F\to F_{C_1}\oplus\hdots\oplus F_{C_N}\to T\to 0$$ where $T$ is a
torsion sheaf whose support is precisely the set of those points
$P_i$, $i=1,\hdots,k$, where $F$ is locally free. From this exact
sequence, we get
\begin{align}
\rk_H(F)&=\frac{1}{h}\sum_{i=1}^{N}
h_{C_i}\rk_{H_{C_i}}(F_{C_i})\notag\\ \dg_H(F)&=\sum_{i=1}^{N}
(\dg_{H_{C_i}}(F_{C_i})+\rk_{H_{C_i}}(F_{C_i})\chi(\mathcal{O}_{C_i}))-\rk_{H}(F)\chi(\mathcal{O}_X)-\chi(T).\notag
\end{align}
In particular, if $F$ is a torsion free sheaf of rank 1 with
respect to $H$, then for every proper subcurve $D$ of $X$ $F_D$ is
torsion free of rank 1 with respect to $H_D$.

 The following lemma, also due to Seshadri, describes the stalk of a
torsion free sheaf on $X$ at the points $P_i$.

\begin{lemma}\label{lem3}
Let $F$ be a pure dimension one sheaf on $X$. If $P_i$ is an
ordinary double point lying in two irreducible components $C_i^1$
and $C_i^2$, then $$F_{P_i}\simeq \mathcal{O}_{X, P_i}^{a_1}\oplus
\mathcal{O}_{C_i^1, P_i}^{a_2}\oplus \mathcal{O}_{C_i^2,
P_i}^{a_3}$$ where $a_1,\ a_2,\ a_3$ are the integer numbers
determined by:
\begin{align}
&a_1+a_2=\rg(F_{P_i}\underset{\mathcal{O}_{X, P_i}}\otimes
\mathcal{O}_{C^1_i,P_i})\notag\\
&a_1+a_3=\rg(F_{P_i}\underset{\mathcal{O}_{X, P_i}}\otimes
\mathcal{O}_{C^2_i,P_i})\notag\\ &a_1+a_2+a_3=\rg(F_{P_i}\otimes
\kappa)\notag
\end{align}
\end{lemma}

\begin{proof} See \cite{Ses}, Huiti\^{e}me Partie, Prop. 3.
\end{proof}

According to the general theory, for every semistable sheaf $F$
with respect to $H$ there is a Jordan-H\"{o}lder filtration
$$0=F_0\subset F_1\subset \hdots \subset F_n=F$$ with stable
quotients $F_i/F_{i-1}$ and $\mu_H(F_i/F_{i-1})=\mu_H(F)$ for
$i=1,\hdots,n$. This filtration need not be unique, but the graded
object $Gr(F)=\textstyle{\bigoplus_{i}} F_i/F_{i-1}$ does not
depend on the choice of the Jordan-H\"{o}lder filtration. Two
semistable sheaves $F$ and $F'$ on $X$ are said to be
$S$-equivalent if $Gr(F)\simeq Gr(F')$. Note that two stable
sheaves are $S$-equivalent only if they are isomorphic.

In the relative case, given a scheme $S$ of finite type over
$\kappa$, a projective morphism of schemes $f\colon X\to S$ whose
geometric fibers are curves and a relative polarization $H$, we
define the relative rank and degree of a coherent sheaf $F$ on
$X$, flat over $S$, as its rank and degree on fibers, and we say
that $F$ is relatively pure of dimension one (resp. stable, resp.
semistable) if it is flat over $S$ and if its restriction to every
geometric fiber of $f$ is pure of dimension one (resp. stable,
resp. semistable).

Let $\bf{Jac^d(X/S)}_s$ (resp. $\bf{\overline{Jac}^d(X/S)}_s$) be
the functor which to any $S$-scheme $T$ associates the set of
equivalence classes of stable invertible (resp. relatively torsion
free rank 1) sheaves on $X_T=X\underset{S}\times T$ with relative
degree $d$. Two such sheaves $F$ and $F'$ are said to be
equivalent if $F'\simeq F\otimes f_T^*N$, where $N$ is a line
bundle on $T$ and $f_T\colon X_T\to T$ is the natural projection.
Similarly, we define the functor $\bf{Jac^d(X/S)}$ (resp.
$\bf{\overline{Jac}^d(X/S)}$) of semistable invertible (resp.
relatively torsion free rank 1) sheaves.

 As a particular case of the Simpson's work
\cite{S}, there exists a projective scheme $\Ja^d(X/S)\to S$ which
coarsely represents the functor $\bf{\overline{Jac}^d(X/S)}$.
Rational points of $\Ja^d(X/S)$ correspond to $S$-equivalence
classes of semistable torsion free sheaves of rank 1 and degree
$d$ on a fiber $X_s$ ($s\in S$). Moreover, $\bf{Jac^d(X/S)}$ is
coarsely represented  by a subscheme $\jac^d(X/S)$ and there are
open subschemes $\ja^d(X/S)_s$ and $\Ja^d(X/S)_s$ which represent
the other two funtors.

\begin{definition} The Simpson Jacobian of $X$ is
$\jac^d(X)_s=\ja^d(X/\Spe \kappa)_s$. We denote
$\Jac^d(X)=\Ja^d(X/\Spe \kappa)$ and $\Jac^d(X)_s=\Ja^d(X/\Spe
\kappa)_s$.
\end{definition}

When $X$ is an integral curve every torsion free rank 1 sheaf is
stable, and then $\jac^d(X)=\jac^d(X)_s$ is equal to the Picard
scheme $\Pic^d(X)$ and $\Jac^d(X)=\Jac^d(X)_s$ coincides with the
Altman-Kleiman's compactification \cite{AK}.

\begin{definition} A generalized tree-like curve is a projective, reduced and connected
curve $X=C_1\cup\hdots\cup C_N$ over $\kappa$ such that the
intersection points, $P_1,\hdots,P_k$, of its irreducible
components are disconnecting ordinary double points.
\end{definition}

Henceforth we shall assume that $X$ is a generalized tree-like
curve, and then  $k=N-1$.

\begin{lemma}\label{lem4}  It is possible to
order the irreducible components $C_1,\hdots,C_N$ of the curve
$X$, so that for every $i\leq N-1$ all but one of the connected
components of $X-C_i$ consists entirely of irreducible components
with subindices smaller than $i$. Then there are irreducible
components, say $C_{i_1},\hdots,C_{i_k}$, with all subindices
smaller than $i$, such that $X_i=C_i\cup C_{i_1}\cup\hdots \cup
C_{i_k}$ is connected and intersects its complement
$\overline{X}_i$ in $X$ in just one point $P_i$.
\end{lemma}

\begin{proof} Teixidor proves this lemma in (\cite{T}, Lem. 1) when $X$ is a
tree-like curve, that is, its irreducible components are smooth,
but her proof is valid for our curve.
\end{proof}

\section{The Description}
Let us suppose from now on that an ordering of the components of
$X$ as in lemma \ref{lem4}, has been fixed.

{\bf Notation:} If $g=\gen(X)$ denote the arithmetic genus of $X$,
that is, the dimension of $\coh^1(X,\mathcal{O}_X)$, for any
torsion free rank 1 sheaf $F$ en $X$ of degree $d$, let $b$ be the
residue class of $d-g$ modulo $h$ so that $$d-g=ht+b.$$ For every
proper subcurve $D$ of $X$, we shall write
$$k_D=\frac{h_D(b+1)}{h}.$$

Lemma \ref{lem4} allows us to define inductively integer numbers
$d_i^X$ as follows:
\begin{align}\label{des1}
 &d_i^X=-\chi(\mathcal{O}_{X_i})+h_{X_i}t+[k_{X_i}]+1-d_{i_1}^X-\hdots-d_{i_k}^X,\
\text{ for }  i=1,\hdots, N-1 \notag\\
&d_N^X=d-d_1^X-\hdots-d_{N-1}^X.
\end{align}

We are now going to modify the above numbers to obtain new numbers
$d_i$ associated with $X$. This is accomplished by a recurrent
algorithm. In order to describe it we start by saying that a
connected subcurve $D=C_{j_1}\cup\hdots\cup C_{j_r}$, $r\geq 1$,
of $X$ ordered according to lemma \ref{lem4} is {\it final} either
when the numbers $k_{D_{j_t}}$ are not integers for $t=1,\hdots,
r-1$ or $D$ is irreducible.

If $D$ is a final curve, we define $d_{j_t}$ as follows:

1. if $r>1$, $d_{j_t}=d^D_{j_t}$ for $t=1,\hdots,r$.

2. if $r=1$,
$d_{j_1}=h_{C_{j_1}}t+[k_{C_{j_1}}]-\chi(\mathcal{O}_{C_{j_1}})$.

\noindent {\bf Algorithm:} If the curve $X$ is final, $d_i=d_i^X$
for all $i$. Otherwise, let $i$ be the first index for which
$k_{X_i}\in \mathbb{Z}$. We consider the two connected components,
$Y=X_i$ and $\overline{Y}=\overline{X_i}$, of $X-P_i$ and we
reorder them according with lemma \ref{lem4}. This induces a new
ordering $P^Y_r$ (resp. $P^{\overline Y}_s$) of the points $P_j$
($j\ne i$) in $Y$ (resp. $\overline{Y}$). Then,

a) If $Y$ (resp. $\overline{Y}$) is a final curve, the process
finishes for $Y$ (resp. $\overline{Y}$).

b) If $Y$ is not final, we take the first index $r$ of $Y$ for
which $k_{Y_r}\in \mathbb{Z}$, consider the connected components,
$Z$ and $\overline{Z}$, of $Y-P^Y_r$ and reorder them according
with lemma \ref{lem4}. If $Z$ and $\overline{Z}$ are final, the
process finishes for $Y$. Otherwise, we iterate the above argument
for those components that are not final and so on. The process
finishes for $Y$ when all subcurves that we find are final.

c) If $\overline{Y}$ is not final, we take the first index $s$ of
$\overline{Y}$ such that $k_{\overline{Y}_s}\in \mathbb{Z}$,
consider the connected components, $W$ and $\overline{W}$, of
$\overline{Y}-P^{\overline{Y}}_s$ and reorder them according with
lemma \ref{lem4}. If $W$ and $\overline{W}$ are final, the process
finishes for $\overline{Y}$. Otherwise, we repeat the above
argument for those components that are not final and so on. The
process finishes for $\overline{Y}$ when all subcurves that we
obtain are final.

The algorithm for $X$ finishes when it finishes for both $Y$ and
$\overline{Y}$. \fin

We can now state the theorem that determines the structure of the
Simpson Jacobian of $X$ and of the schemes $\Jac^d(X)_s$ and
$\Jac^d(X)$.

\begin{theorem} \label{t} Let $X=C_1\cup\hdots\cup C_N$ , $N\geq 2$, be a
generalized tree-like curve.

a) If $k_{X_i}$ is not an integer for every $i\leq N-1$, then
$$\jac^d(X)_s=\prod_{i=1}^{N}\Pic^{d_i^X}(C_i)\ \ \ \text{and}$$
$$\jac^d(X)_s\subseteq \Jac^d(X)_s=\Jac^d(X)\simeq \prod_{i=1}^{N}
\Ja^{d_i^X}(C_i)$$ where $d_i^X$ are the above integers.

b) If $k_{X_i}$ is an integer for some $i\leq N-1$, then
$$\jac^d(X)_s=\Jac^d(X)_s=\emptyset\ \ \ \text{and}$$
$$\Jac^d(X)\simeq \prod_{i=1}^{N}\Ja^{d_i}(C_i)$$ where $d_i$ are
the integers constructed with the above algorithm.

\end{theorem}

In order to prove that the Simpson Jacobian of $X$ is not empty
only when $k_{X_i}\notin \mathbb{Z}$ for all $i\leq N-1$, we need
the next two lemmas that characterize stable invertible sheaves on
$X$.

\begin{lemma}\label{lem5} Let $L$ be a line bundle on $X$ of degree
$d$. If $L$ is stable, then $k_{X_i}$ is not an integer for every
$i\leq N-1$ and $L$ is obtained by gluing invertible sheaves $L_i$
on $C_i$ of degrees $d_i^X$, $i=1,\hdots,N$.
\end{lemma}

\begin{proof} Let us consider the subsheaves $L^{X_i}$ of $L$,
$i=1,\hdots,N-1$, $X_i$ being the subcurves of $X$ given by lemma
\ref{lem4}. By the stability of $L$, we get
\begin{equation}\label{des2}
h_{X_i}d-hd_{X_i}\des
h\chi(\mathcal{O}_{X_i})-h_{X_i}\chi(\mathcal{O}_X)\ \ \text{ for
} i=1,\hdots, N-1 \tag{2}
\end{equation}
Considering the subsheaves  $L^{\overline{X}_i}$ of $L$,
$i=1,\hdots,N-1$, yields
\begin{equation}\label{des3}
h_{\overline{X}_i}d-hd_{\overline{X}_i}\des
h\chi(\mathcal{O}_{\overline{X}_i})-h_{\overline{X}_i}\chi(\mathcal{O}_X)\
\ \text{ for } i=1,\hdots, N-1 \tag{3}
\end{equation}
Since $X=X_i\cup \overline{X}_i$ and $X_i$, $\overline{X_i}$ meet
only at $P_i$, we have $d=d_{X_i}+d_{\overline{X}_i}$,
$h=h_{X_i}+h_{\overline{X}_i}$ and
$\chi(\mathcal{O}_X)=\chi(\mathcal{O}_{X_i})+\chi(\mathcal{O}_{\overline{X}_i})-1$.
Then, \eqref{des2} and \eqref{des3} give
\begin{equation}\label{des4}
h_{X_i}d+h_{X_i}\chi(\mathcal{O}_X)-h\chi(\mathcal{O}_{X_i})\des
hd_{X_i}\des
h_{X_i}d+h_{X_i}\chi(\mathcal{O}_X)-h\chi(\mathcal{O}_{X_i})+h
\tag{4}
\end{equation}
We have $\rk_{H_{X_i}}(L_{X_i})=1$ so that
$d_{X_i}=\dg_{H_{X_i}}(L_{X_i})$ is an integer. Then, if
$k_{X_i}\in \mathbb{Z}$ for some $i\leq N-1$, \eqref{des4} becomes
a contradiction. Thus $k_{X_i}\notin \mathbb{Z}$ for all $i\leq
N-1$ and there is only one possibility for $d_{X_i}$, namely
$$d_{X_i}=-\chi(\mathcal{O}_{X_i})+h_{X_i}t+[k_{X_i}]+1,\ \ \
\text{ for } i=1,\hdots, N-1. $$ From
$d_{X_i}=d_{C_i}+d_{C_{i_1}}+\hdots+d_{C_{i_k}}$ and the exact
sequence $$0\to L\to L_{C_1}\oplus\hdots\oplus L_{C_N}\to
\oplus_{i=1}^{N-1}\kappa(P_i)\to 0$$ we deduce that
$d_{C_i}=d_i^X$ for all $i$ and the proof is complete.
\end{proof}

The following lemma, which is similar to lemma 2 of Teixidor
\cite{T}, proves the converse of lemma \ref{lem5}.

\begin{lemma}\label{lem6} Let $L$ be an invertible sheaf of degree
$d=g+ht+b$ obtained by gluing line bundles $L_i$ on $C_i$ of
degrees $d_i^X$, $i=1,\hdots,N$. Suppose that $k_{X_i}$ is not
integer for every $i\leq N-1$ and let $D=C_{a(1)}\cup\hdots\cup
C_{a(t)}$ be a proper subcurve of $X$. Then,

1. The following inequality holds:
$$-\chi(\mathcal{O}_D)+h_Dt+k_D< d_D<
-\chi(\mathcal{O}_D)+h_Dt+k_D+\alpha$$ where $\alpha$ is the
number of intersection points of $D$ and its complement $\overline
D$ in $X$.

2. $L$ is stable.
\end{lemma}

\begin{proof} 1. If $D$ is a connected subcurve of $X$, let us denote by
$C_{i(k)}$, $k=1,\hdots, \alpha$, (resp. $C_{j(k)}$) the component
of $D$ (resp. $\overline{D}$) which contains the point $P_k$ (it
is possible to have $C_{i(k)}=C_{i(k')}$ for $k \ne k'$).

Suppose first that $i(k)>j(k)$ for $k=1,\hdots, \alpha$. Then, the
subcurve $X_{j(k)}$ of lemma \ref{lem4} is the connected component
of $X-P_k$ which contains $C_{j(k)}$.

Let us see that $\overline{D}=\bigsqcup_{k=1}^{\alpha}X_{j(k)}$.
Since the inclusion
$\overline{D}\subseteq\bigcup_{k=1}^{\alpha}X_{j(k)}$ is clear we
have only to prove that $\bigsqcup_{k=1}^{\alpha}X_{j(k)}$
contains no component of $D$ and that $X_{j(k)}\cap X_{j(k')}=
\emptyset$ if $k\ne k'$ . If $C_{i(k')}\subseteq X_{j(k)}$ for
some $k'$, then $C_{j(k')}\subseteq X_{j(k)}$ and
$j(k')<i(k')<j(k)<i(k)$. Hence, $$X-P_k-P_{k'}=(Z_k^1\sqcup
Z_k^2)\sqcup \overline{X_{j(k)}}$$ where we denote by $Z_k^1$
(resp. $Z_k^2$) the connected component of $X_{j(k)}-P_{k'}$ that
contains $C_{i(k')}$ (resp. $C_{j(k')}$). Analogously, $C_{j(k)}$,
$C_{i(k)}$ are contained in $\overline{X_{j(k')}}$ and
$$X-P_{k'}-P_k=X_{j(k')}\sqcup (Z_{k'}^1\sqcup Z_{k'}^2)$$ where
$Z_{k'}^1$ (resp. $Z_{k'}^2$) is the connected component of
$\overline{X_{j(k')}}-P_k$ containing $C_{i(k)}$ (resp.
$C_{j(k)}$). We deduce that $X_{j(k')}=Z^2_k$,
$\overline{X_{j(k)}}=Z^1_{k'}$ and $Z^2_{k'}=Z^1_k$. Thus,
$C_{j(k)}\subseteq Z^1_k$ and $j(k)<i(k')$, which is absurd.
Therefore, none of the components $C_{i(k')}$ of $D$ is in
$\bigcup_{k=1}^{\alpha}X_{j(k)}$ and, since $X_{j(k)}$ are
connected, no other component of $D$ is either. Moreover, a
similar argument shows that $X_{j(k)}$ and $X_{j(k')}$ have no
common components for $k\ne k'$ .

Let now $t_k+1$ be the number of irreducible components in
$X_{j(k)}$. Then, $\sum_{k} (t_k+1)$ is equal to $N-t$. On the
other hand, the number of intersection points in
$\bigsqcup_{k}X_{j(k)}$ is $\sum_{k} t_k+s$, where $s$ is the
number of intersection points of $X_{j(k)}$, $k=1,\hdots, \alpha$.
Since $\sum_{k} t_k+s$ is the number $N-t-\alpha$ of intersection
points in $\overline{D}$, we get $s=0$ and
$\overline{D}=\bigsqcup_{k=1}^{\alpha}X_{j(k)}$. Then,
$d_{\overline D}=\sum_k d_{X_{j(k)}}$, $h_{\overline D}=\sum_k
h_{X_{j(k)}}$ and
$\chi(\mathcal{O}_{\overline{D}})=\sum_k\chi(\mathcal{O}_{X_{j(k)}}).$

By definition of $d_i^X$, for every $k=1,\hdots,\alpha$, we have
\begin{equation}\label{des5}
-\chi(\mathcal{O}_{X_{j(k)}})+h_{X_{j(k)}}t+k_{X_{j(k)}}<d_{X_{j(k)}}<
-\chi(\mathcal{O}_{X_{j(k)}})+h_{X_{j(k)}}t+k_{X_{j(k)}}+1\tag{5}
\end{equation}
Taking into account that $d=d_D+d_{\overline D}$,
$h=h_D+h_{\overline D}$ and
$\chi(\mathcal{O}_X)=\chi(\mathcal{O}_{D})+\chi(\mathcal{O}_{\overline
D})-\alpha$, we obtain the result in this case.

Assume now that $j(1)>i(1)$. Then, $X_{i(1)}$ is the connected
component of $X-P_1$ which contains $C_{i(1)}$. Since $D$ is
connected, $X_{i(1)}$ contains $D$  and then $X_{i(1)}$ contains
$C_{i(k)}$ and $C_{j(k)}$ for $k=2,\hdots, \alpha$. This implies
that $i(k)<i(1)$ for $k=2,\hdots, \alpha$. If we had $i(k)<j(k)$
for some $k=2,\hdots, \alpha$, then $D$ would be contained in
$X_{i(k)}$ as before and $i(1)<i(k)$, which is absurd. Thus,
$i(k)>j(k)$ for $k=2,\hdots, \alpha$. In this situation, one has
that $X_{i(1)}=D\cup Z$, where $Z=\bigcup_{k=2}^\alpha X_{j(k)}$
and $D$ intersects $Z$ in $\alpha-1$ points.

Thus, we have relations \eqref{des5} for $k=2,\hdots, \alpha$ and,
arguing as in the former case, we obtain
\begin{equation}\label{des6}
-\chi(\mathcal{O}_Z)+h_Zt+k_Z< d_Z<
-\chi(\mathcal{O}_Z)+h_Zt+k_Z+(\alpha-1)\tag{6}
\end{equation}
By definition of $d_i^X$, one has
\begin{equation}
-\chi(\mathcal{O}_{X_{i(1)}})+h_{X_{i(1)}}t+k_{X_{i(1)}}<d_{X_{i(1)}}<
-\chi(\mathcal{O}_{X_{i(1)}})+h_{X_{i(1)}}t+k_{X_{i(1)}}+1\notag
\end{equation}
which together with \eqref{des6} proves the statement.

If the subcurve $D$ is not connected, the inequality holds for
every connected component and then it is easy to deduce it for
$D$.

2. By lemma \ref{lem2}, it is enough to show that $\mu_{H}(L^D)<d$
for every proper subcurve $D$ of $X$. Since
$$\mu_H(L^D)=\frac{hd-hd_D+h_D\chi(\mathcal{O}_X)-h\chi(\mathcal{O}_D)}{h-h_D},$$
the result follows from 1.
\end{proof}

\begin{proof}[Proof of Theorem~\ref{t}] By lemmas \ref{lem5} and \ref{lem6}, $\jac^d(X)_s$ is
not empty only if $k_{X_i}$ is not integer for every $i\leq N-1$,
and in this case it is equal to $\prod_{i=1}^N\Pic^{d_i^X}(C_i)$.

We prove now the remaining statements of the theorem. If $L$ is a
strictly semistable line bundle on $X$ of degree $d$ then,
$hd_{X_i}$ is equal to one of the two extremal values of the
inequality \eqref{des4}. In particular, $k_{X_i}$ is an integer
for some $i\leq N-1$.

Let $i$ be the first index such that $k_{X_i}$ is integer. Then,
there are two possibilities for $d_{X_i}$:

a) $d_{X_i}=-\chi(\mathcal{O}_{X_i})+h_{X_i}t+k_{X_i}$

b) $d_{X_i}=-\chi(\mathcal{O}_{X_i})+h_{X_i}t+k_{X_i}+1$

Let us construct a Jordan-H\"{o}lder filtration for $L$ in both cases.
Since case a) and case b) are the same but with the roles of $X_i$
and $\overline{X_i}$ intertwined, we give the construction in the
case a).

We have that $\mu_H(L_{X_i})=\mu_H(L^{X_i})=\mu_H(L)$. Then,
$L_{X_i}$ and $L^{X_i}\simeq L_{\overline{X_i}}(-P_i)$ are
semistable with respect to $H$ and, by lemma \ref{lem1}, they are
semistable with respect to $H_{X_i}$ and $H_{\overline{X_i}}$
respectively.

For simplicity, we shall write  $Y=X_i=C_{i_0}\cup
C_{i_1}\cup\hdots \cup C_{i_k}$ with $i_1,\hdots, i_k<i_0=i$ and
$Z=\overline{X_i}$, which are again generalized tree-like curves.

Let us see when the sheaves $L_Y$ and $L_Z(-P_i)$ are stable. We
can fix an ordering for $Y$ as in lemma \ref{lem4}, so that
$Y=C_{\sigma(i_0)}\cup\hdots\cup C_{\sigma(i_k)}$ and we obtain
subcurves  $Y_r$ of $Y$ for
$r=\sigma(i_0),\hdots,\sigma(i_{k-1})$.

{\it Claim 1.} The sheaf $L_Y$ is stable if and only if $k_{Y_r}$
is not an integer for $r=\sigma(i_0),\hdots,\sigma(i_{k-1})$.

\demo Since the residue class of $d_Y-\gen(Y)$ modulo $h_Y$ is
$b_Y=k_Y-1$, the numbers $\frac{h_{Y_r}(b_Y+1)}{h_Y}=k_{Y_r}$ are
not integers for $r=\sigma(i_0),\hdots,\sigma(i_{k-1})$. Then,
from lemma \ref{lem6}, we have only to prove that $L_Y$ is in
$\prod_r \Pic^{d_r^Y}(C_r)$, where $d_r^Y$ are the integer numbers
defined as $d_i^X$ but with the new ordering of $Y$ and $r$ runs
through the irreducible components of $Y$. This is equivalent to
proving that
$$d_{Y_r}=-\chi(\mathcal{O}_{Y_r})+h_{Y_r}t+[k_{Y_r}]+1\ \ \
\text{ for } r=\sigma(i_0),\hdots,\sigma(i_{k-1}).$$ Actually,
since $L$ is semistable and $Y_r$ is a proper subcurve of $X$,
arguing as in lemma \ref{lem6}, we obtain
\begin{equation}\label{des7}
-\chi(\mathcal{O}_{Y_r})+h_{Y_r}t+k_{Y_r}\leq d_{Y_r}\leq
-\chi(\mathcal{O}_{Y_r})+h_{Y_r}t+k_{Y_r}+\alpha\tag{7}
\end{equation}
where $\alpha$ is the number of intersection points of $Y_r$ and
its complement in $X$. We have that $\alpha\leq 2$ and $d_{Y_r}$
is not equal to the extremal values of \eqref{des7} because
$k_{Y_r}\notin \mathbb{Z}$. Moreover, if it were
$$d_{Y_r}=-\chi(\mathcal{O}_{Y_r})+h_{Y_r}t+[k_{Y_r}]+2,$$ since
$d_Y=d_{Y_r}+d_{\overline{Y_r}^Y}$,
$h_Y=h_{Y_r}+h_{\overline{Y_r}^Y}$ and
$\chi(\mathcal{O}_Y)=\chi(\mathcal{O}_{Y_r})+\chi(\mathcal{O}_{\overline{Y_r}^Y})-1$,
$\overline{Y_r}^Y$ being the complement of $Y_r$ in $Y$, then
$$d_{\overline{Y_r}^Y}=-\chi(\mathcal{O}_{\overline{Y_r}^Y})+h_{\overline{Y_r}^Y}t+
[k_{\overline{Y_r}^Y}]$$ which contradicts  the semistability of
$L$. Thus,
$$d_{Y_r}=-\chi(\mathcal{O}_{Y_r})+h_{Y_r}t+[k_{Y_r}]+1$$ and the
proof of the claim 1 is complete.

On the other hand, the irreducible components of $Z$ are ordered
according the instructions in lemma \ref{lem4} and the subcurves
$Z_s$, where $s$ runs through the irreducible components of $Z$
and $s\leq N-1$, are equal to either $X_s$ or $X_s-Y$.

{\it Claim 2.} The sheaf $L_Z(-P_i)$ is stable if and only if
$k_{X_s}$ is not an integer for every $s>i$.

\demo Since the residue class of $d_Z-1-\gen(Z)$ modulo $h_Z$ is
$b_Z=k_Z-1$ and $k_Y\in \mathbb{Z}$ , by the hypothesis, the
numbers $\frac{h_{Z_s}(b_Z+1)}{h_Z}$ are not integers for $s>i $
and, by the choice of $i$, they aren't for $s<i$ either. Then, by
lemma \ref{lem6}, it is enough to prove that
$$d_{H_{Z_s}}(L_Z(-P_i)|_{Z_s})=-\chi(\mathcal{O}_{Z_s})+h_{Z_s}t+[k_{Z_s}]+1\
\ \ \text{ for } s\leq N-1.$$ Since $L$ is semistable, we have
that $$d_{X_s}=-\chi(\mathcal{O}_{X_s})+h_{X_s}t+[k_{X_s}]+1.$$
Moreover, if  $Z_s=X_s$ then,
$d_{H_{Z_s}}(L_Z(-P_i)|_{Z_s})=d_{X_s}$ and if $Z_s=X_s-Y$ then,
$d_{H_{Z_s}}(L_Z(-P_i)|{Z_s})=d_{X_s}-d_Y-1$. We obtain the
desired result in both cases and the proof of the claim 2 is
complete.

We return now to the proof of the theorem. If $k_{Y_r}$ and
$k_{X_s}$ are not integers for
$r=\sigma(i_0),\hdots,\sigma(i_{k-1})$ and $s>i$, then $0\subset
L_Z(-P_i)\subset L$  is a Jordan-H\"{o}lder filtration for $L$ and the
$S$-equivalence class of $L$ belongs to
$\prod_r\Pic^{d_r^Y}(C_r)\times \prod_s\Pic^{d_s^X}(C_s)$.

On the other hand, if $k_{Y_r}$ is integer for some
$r=\sigma(i_0),\hdots,\sigma(i_{k-1})$, the sheaf $L_Y$ is
strictly semistable and we have to repeat the above procedure with
$L_Y$ in the place of $L$ and the curve $Y$ in the place of $X$.
Similarly, if $k_{X_s}$ is integer for some $s>i$, the sheaf
$L_Z(-P_i)$ is strictly semistable. Then, we have to repeat the
above procedure for $L_Z(-P_i)$. By iterating this procedure, we
get a Jordan-H\"{o}lder filtration for $L_Y$: $$0=F_0\subset
F_1\subset \hdots \subset F_m=L_Y$$ and another for $L_Z(-P_i)$:
$$0=G_0\subset G_1\subset \hdots\subset G_n=L_Z(-P_i).$$
Therefore, a filtration for $L$ is given by $$0=G_0\subset
G_1\subset \hdots\subset L_Z(-P_i)\subset \pi_Y^{-1}(F_1)\subset
\hdots\subset \pi_Y^{-1}(L_Y)=L.$$ Thus, the $S$-equivalence class
of $L$ belongs to $\prod_{i=1}^{N}\Pic^{d_i}(C_i)$, where $d_i$
are the integer numbers constructed with the algorithm.

Finally, let us consider a torsion free sheaf $F$ on $X$ of rank 1
and degree $d$ which is not locally free. When $F$ is locally free
at the intersection points $P_i$ for all $i=1,\hdots,N-1$,
calculations and results are analogous to the former ones. If $F$
is not locally free at $P_i$ for some $i=1,\hdots,N-1$, then there
is a natural morphism $$F\to F_Y\oplus F_Z$$ where $Y$, $Z$ are
the connected components of $X-{P_i}$, that is clearly an
isomorphism outside $P_i$. But this is an isomorphism at $P_i$ as
well because by lemma \ref{lem3}, $F_{P_i}\simeq
\mathcal{O}_{C_i^1, P_i}\oplus \mathcal{O}_{C_i^2, P_i}$ and this
is precisely the stalk of $F_Y\oplus F_Z$ at $P_i$. We conclude
that if $F$ is strictly semistable, then $k_Y$ y $k_Z$ are
integers, $d_Y$ and $d_Z$ are given by
$$d_Y=-\chi(\mathcal{O}_Y)+h_Yt+k_Y,\quad\quad
d_Z=-\chi(\mathcal{O}_Z)+h_Zt+k_Z,$$ and $F_Y$ and $F_Z$  are
semistable as well. Then, the construction of  a Jordan-H\"{o}lder
filtration for $F$ can be done as above and thus the
$S$-equivalence class of $F$ belongs to
$\prod_{i=1}^{N}\Jac^{d_i}(C_i)$.
\end{proof}

We now give three examples to illustrate this theorem.

{\bf Example 1.} Let $X=C_1\cup\hdots\cup C_N$, $N\geq 2$, be a
generalized tree-like curve with a polarization $H$ whose degree
$h$ is a prime number. Suppose that the irreducible components of
$X$ are ordered according with lemma \ref{lem4}. Then, since
$h_{X_i}$ is not divisible by $h$,
$k_{X_i}=\frac{h_{X_i}(b+1)}{h}$ is an integer if and only if
$b=h-1$. Therefore, by theorem \ref{t}, we have that if $b<h-1$,
$\jac^d(X)_s=\prod_{i=1}^{N}\Pic^{d_i^X}(C_i)$ and
$$\jac^d(X)_s\subsetneqq \Jac^d(X)_s=\Jac^d(X)\simeq
\prod_{i=1}^{N} \Ja^{d_i^X}(C_i),$$ whereas for $b=h-1$,
$\jac^d(X)_s$ and $\Jac^d(X)_s$ are empty and $$\Jac^d(X)\simeq
\prod_{i=1}^{N}\Ja^{d_i}(C_i).$$ Moreover, in this case the number
$d_i$ is given by $$d_i=h_{C_i}(t+1)-\chi(\mathcal{O}_{C_i}) \ \ \
\text{ for } i=1,\hdots,N.$$ Actually, if $L$ is a strictly
semistable line bundle of degree $d$ on $X$, the first index $i$
such that $k_{X_i}$ is integer is $i=1$ and the connected
components of $X-P_1$ are $Y=C_1$ and $Z=C_2\cup \hdots \cup C_N$.
Suppose, as in the proof of the theorem, that
$$d_Y=-\chi(\mathcal{O}_Y)+h_Yt+k_Y=h_Y(t+1)-\chi(\mathcal{O}_Y).$$
Then, $L_Y$ and $L_Z(-P_1)$ are semistable and $L_Y=L_{C_1}\in
\Pic^{h_{C_1}(t+1)-\chi(\mathcal{O}_{C_1})}(C_1)$. On the other
hand, since $k_{X_s}$ is integer for $s=2,\hdots,N$, we have to
apply the procedure to the curve $Z$. Here, the first index $s$
such that $k_{X_s}$ is integer is $s=2$ and the connected
components of $Z-P_2$ are $C_2$ and $C_3\cup\hdots\cup C_N$.
Proceeding exactly as in case of $X$, we obtain a sheaf belonging
to $\Pic^{h_{C_2}(t+1)-\chi(\mathcal{O}_{C_2})}(C_2)$ and another
supported on $C_3\cup\hdots\cup C_N$. Since $k_{X_s}$ is integer
for $s=3,\hdots,N$, we apply the procedure to the curve
$C_3\cup\hdots\cup C_N$, and so on. The iteration of this
procedure will only finish when we obtain a sheaf supported on
$C_N$ that belongs to
$\Pic^{h_{C_N}(t+1)-\chi(\mathcal{O}_{C_N})}(C_N)$. Thus,
$d_i=h_{C_i}(t+1)-\chi(\mathcal{O}_{C_i})$ for $i=1,\hdots,N$.

{\bf Example 2.} We are now going to recover examples 1 and 2 of
\cite{A}. There, the irreducible components $C_i$ are taken to be
smooth and $d=g-1$. Then, the residue class of $d-g$ modulo $h$ is
$b=h-1$ and $t=-1$. It follows that $k_{X_i}$ is integer for
$i=1,\hdots,N$ and, arguing as in our example 1, we have that
$\jac^{g-1}(X)_s$ and $\Jac^{g-1}(X)_s$ are empty and
$$\Jac^{g-1}(X)\simeq
\prod_{i=1}^N\Jac^{h_{C_i}(t+1)-\chi(\mathcal{O}_{C_i})}(C_i)\simeq
\prod_{i=1}^N\Pic^{g_i-1}(C_i)$$ as asserted in \cite{A}.

{\bf Example 3.} Let $X$ be the following tree-like curve

\vspace{1truecm}
\hspace{4truecm}\includegraphics[scale=0.5]{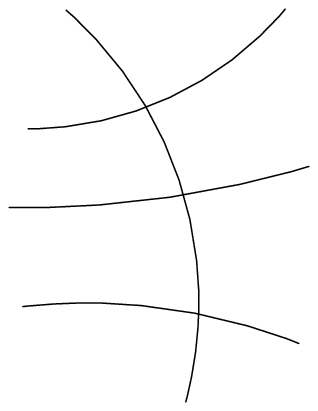}
\vspace{1truecm}

Fixing an ordering of the irreducible components of $X$ as in
lemma \ref{lem4}, we obtain

\vspace{1truecm}
\hspace{4truecm}\includegraphics[scale=0.5]{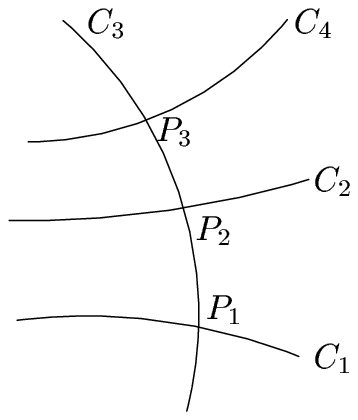}
\vspace{1truecm}

\noindent so that $X_1=C_1$, $X_2=C_2$, $X_3=C_1\cup C_2\cup C_3$.
Assume that the first index $i$ such that $k_{X_i}$ is integer is
$i=3$ and let us compute the numbers $d_i$, $i=1,\hdots,4$ in this
case.

The connected components of $X-P_3$ are $Y=C_1\cup C_2\cup C_3$
and $Z=C_4$. Suppose $L$ is a strictly semistable line bundle of
degree $d$ on $X$ and that $$d_Y=-\chi(\mathcal{O}_Y)+h_Yt+k_Y$$
(the other case is similar). Then $L_Y$ and $L_Z(-P_3)$ are
semistable. Moreover, $L_Z(-P_3)\in
\Pic^{h_{C_4}t+k_{C_4}-\chi(\mathcal{O}_{C_4})}(C_4)$ so that
$$d_4=h_{C_4}t+k_{C_4}-\chi(\mathcal{O}_{C_4}).$$  We now have to
fix a new ordering for $Y$ according with lemma \ref{lem4}. We can
take, for instance, $Y=C_{\sigma(1)}\cup C_{\sigma(2)}\cup
C_{\sigma(3)}$ with $\sigma(1)=1$, $\sigma(2)=3$ and
$\sigma(3)=2$. Then, $Y_{\sigma(1)}=C_1$ and
$Y_{\sigma(2)}=C_1\cup C_3$. Therefore, since
$k_{Y_{\sigma(1)}}=k_{X_1}$ and $k_{Y_{\sigma(2)}}=k_Y-k_{X_2}$
are not integers, we conclude that
\begin{align}
d_1&=d^Y_{\sigma(1)}=-\chi(\mathcal{O}_{Y_{\sigma(1)}})+h_{Y_{\sigma(1)}}t+
[k_{Y_{\sigma(1)}}]+1=d^X_1,\notag\\
d_3&=d_{\sigma(2)}^Y=-\chi(\mathcal{O}_{Y_{\sigma(2)}})+h_{Y_{\sigma(2)}}t+
[k_{Y_{\sigma(2)}}]+1-d^Y_{\sigma(1)}=\notag\\
&=-\chi(\mathcal{O}_Y)+h_Yt+k_Y+1-d^X_1-d^X_2-1=d^X_3-1,\notag\\
d_2&=d_{\sigma(3)}^Y=d_Y-d_{\sigma(1)}^Y-d_{\sigma(2)}^Y=d^X_2.\notag
\end{align}

\end{document}